\documentclass[reqno,english]{smfart}
\usepackage[english]{babel}
\usepackage{smfthm}
\NumberTheoremsIn{subsection}

\author{Yves Laurent}
\address{ Institut Fourier Math\'ematiques \\
          UMR 5582  CNRS/UJF\\
                 BP 74\\
         38402 St Martin d'H\`eres Cedex}
\email{Yves.Laurent@ujf-grenoble.fr}
\urladdr{http://www-fourier.ujf-grenoble.fr/~laurenty}
\title{Regularity of $\mathcal D$-modules associated to a symmetric pair}
\alttitle{R\'egularit\'e des $\mathcal D$-modules associ\'es \`a une paire sym\'etrique}

\DeclareMathOperator{\HOM}{\mathcal Hom}

\begin{document}
\frontmatter
\begin{abstract}
The invariant eigendistributions on a reductive Lie algebra are
solutions of a holonomic $\mathcal D$-module which has been proved to be
regular by Kashiwara-Hotta. We solve here a conjecture of Sekiguchi
saying that in the more general case of symmetric pairs, the
corresponding module are still regular.
\end{abstract}

\begin{altabstract}
Sur une alg\`ebre de Lie r\'eductive, les distributions invariantes qui
sont vecteurs propres des op\'erateurs diff\'erentiels bi-invariants sont
les solutions d'un syst\`eme holonome. Il a \'et\'e d\'emontr\'e par
Kashiwara-Hotta que ce module est r\'egulier. Nous r\'esolvons ici une
conjecture de  Sekiguchi en montrant que ce r\'esultat est encore vrai dans le
cas plus g\'en\'eral des paires sym\'etriques.
\end{altabstract}

\altkeywords{D-module, groupe de Lie, paire sym\'etrique, r\'egularit\'e}
\dedicatory{\`A Jean-Michel Bony, \`a l'occasion de son 60e anniversaire.}
\maketitle

\newcommand{\x}{\times}
\newcommand{\ox}{\otimes}
\renewcommand{\tt}{{\tilde t}}  
\newcommand{\ts}{{\gs}}  
\newcommand{\abs}[1]{\lvert#1\rvert}
\newcommand{\norm}[1]{\lVert#1\rVert}
\newcommand{\ensemble}[2]{\{\,#1\mid#2\,\}}

\newcommand{\C}{{\mathbb C}}
\renewcommand{\L}{{\mathbb L}}
\newcommand{\N}{{\mathbb N}}
\renewcommand{\P}{{\mathbb P}}
\newcommand{\Q}{{\mathbb Q}}
\newcommand{\R}{{\mathbb R}}
\newcommand{\Z}{{\mathbb Z}}
\newcommand{\bN}{{\mathbf N}}

\newcommand{\CA}{\mathcal A}
\newcommand{\CB}{\mathcal B}
\newcommand{\CC}{\mathcal C}
\newcommand{\CaD}{\mathcal D}   
\newcommand{\CE}{\mathcal E}       \newcommand{\cE}{\widehat E}
\newcommand{\CF}{\mathcal F}
\newcommand{\CH}{\mathcal H}
\newcommand{\CI}{\mathcal I}
\newcommand{\CJ}{\mathcal J}
\newcommand{\CL}{\mathcal L}
\newcommand{\CM}{\mathcal M}
\newcommand{\CN}{\mathcal N}
\newcommand{\CP}{\mathcal P}
\newcommand{\CO}{\mathcal O}
\newcommand{\CS}{\mathcal S}
\newcommand{\CW}{\mathcal W}

\newcommand{\kC}{\mathfrak C}
\newcommand{\kG}{\mathfrak G}
\newcommand{\kH}{\mathfrak H}
\newcommand{\kN}{\mathfrak N}
\newcommand{\kO}{\mathfrak O}
\newcommand{\kS}{\mathfrak S}

\newcommand{\kb}{\mathfrak b}
\newcommand{\kc}{\mathfrak c}
\newcommand{\kg}{\mathfrak g}
\newcommand{\kh}{\mathfrak h}
\newcommand{\kk}{\mathfrak k}
\newcommand{\kl}{\mathfrak l}
\newcommand{\kn}{\mathfrak n}
\newcommand{\kp}{\mathfrak p}
\newcommand{\kq}{\mathfrak q}
\newcommand{\ks}{\mathfrak s}
\newcommand{\ku}{\mathfrak u}
\newcommand{\kz}{\mathfrak z}

\newcommand{\kgp}{\kg_\gP}
\newcommand{\kqp}{\kq_\gP}
\newcommand{\khp}{\kh_\gP}
\newcommand{\khpo}{\khp^\bot}
\newcommand{\khpp}{(\khpo)'}
\newcommand{\sld}{\ks\kl_2}
\newcommand{\kps}{\kp^*}

\newcommand{\tQ}{\widetilde Q}
\newcommand{\tq}{\widetilde q}
\newcommand{\tgr}{\widetilde \varrho}

\newcommand{\vespa}{\vspace{1em}}
\renewcommand{\setminus}{-}     
\newcommand{\vf}{\varphi}
\newcommand{\dsur}[1]{\frac \partial{\partial#1} }               
\newcommand{\TX}{{T^*X}}
\newcommand{\TY}{{T^*Y}}
\newcommand{\TL}{{T^*\Lambda}}
\newcommand{\tYX}{{T_Y\!X}}
\newcommand{\pYX}{{P_Y\!X}}
\newcommand{\TYX}{{T^*_YX}}
\newcommand{\dTX}{{\dot T^*X}}
\newcommand{\dTYX}{{\dot T^*_YX}}
\newcommand{\OX}{{{\mathcal O_X}}}
\newcommand{\OY}{{\mathcal O_Y}}
\newcommand{\OZ}{{\mathcal O_Z}}
\newcommand{\DX}{{\mathcal D_X}}
\newcommand{\DXi}[1]{{\mathcal D_{X,#1}}}
\newcommand{\EX}{{{\mathcal E_X}}}
\newcommand{\DY}{{{\mathcal D_Y}}}
\newcommand{\DZ}{{{\mathcal D_Z}}}
\newcommand{\DYX}{{\mathcal D_{Y\rightarrow X}}}
\newcommand{\DXY}{{\mathcal D_{X\leftarrow Y}}}
\newcommand{\DXtY}{{\mathcal D_{X\rightarrow Y}}}
\newcommand{\uYX}{{1_{Y\rightarrow X}}}
\newcommand{\DtYX}{{\mathcal D_{\tYX}}}
\newcommand{\DtYXr}{{\mathcal D_{[\tYX/Y]}}}
\newcommand{\DpYX}{{\mathcal D_{\pYX}}}
\newcommand{\DIXY}{{\mathcal D^\infty_{X\leftarrow Y}}}
\newcommand{\DIY}{{{\mathcal D^\infty_Y}}}
\newcommand{\infun}{{\scriptstyle (\infty,1)}}

\newcommand{\tL}{\widetilde L}
\newcommand{\tM}{\widetilde M}
\newcommand{\tR}{\widetilde R}
\newcommand{\tX}{\widetilde X}
\newcommand{\tY}{\widetilde Y}
\newcommand{\tZ}{\widetilde Z}

\newcommand{\CHM}{{Ch(\CM)}}
\newcommand{\CHMF}{{Ch(\CM_F)}}
\newcommand{\TCHM}{{\widetilde{Ch}(\CM)}}

\renewcommand{\projlim}[1]{\underset{#1}{\varprojlim}}
\renewcommand{\injlim}[1]{\underset{#1}{\varinjlim}}

\def\lacute{\mathopen{<}}
\def\racute{\mathopen{>}}
\newcommand{\scal}[2]{{\lacute#1,#2\racute}}

\newcommand{\ga}{\alpha}
\newcommand{\gb}{\beta}
\renewcommand{\gg}{\gamma}           \newcommand{\gG}{\Gamma}
\newcommand{\gd}{\delta}           \newcommand{\gD}{\Delta}
                                   \newcommand{\tgD}{\widetilde\Delta}
\newcommand{\gep}{\varepsilon}
\newcommand{\gz}{\zeta}
\newcommand{\gh}{\eta}
\newcommand{\gk}{\kappa}
\newcommand{\gth}{\theta}          \newcommand{\gTh}{\Theta}
\newcommand{\gvt}{\vartheta}
\newcommand{\gl}{\lambda}          \newcommand{\gL}{\Lambda}
\newcommand{\gm}{\mu}
\newcommand{\gn}{\nu}
\newcommand{\gx}{\xi}              \newcommand{\gX}{\Xi}
\newcommand{\gp}{\pi}              \newcommand{\gP}{\Pi}
\newcommand{\gpb}{\varpi}
\newcommand{\gro}{\varrho}    
\newcommand{\gs}{\sigma}           \newcommand{\gS}{\Sigma}
\newcommand{\gt}{\tau}
\newcommand{\gf}{\varphi}          \newcommand{\gF}{\Phi}
\newcommand{\gq}{\chi}
\newcommand{\gy}{\psi}             \newcommand{\gY}{\Psi}
\newcommand{\go}{\omega}           \newcommand{\gO}{\Omega}

\newcommand{\tgh}{{\widetilde\eta}}

\newcommand{\exacte}[3]{0\xrightarrow{\ \ }{#1}\xrightarrow{\ \ }{#2}
 \xrightarrow{\ \ }{#3} \xrightarrow{\ \ }0}

\newcommand{\tg}[1]{<\![#1,x],\partial_x\!>}
\newcommand{\tgA}{{\tg A}}
\newcommand{\Dg}{\CaD_\kg}
\newcommand{\Dh}{\CaD_\kh}
\newcommand{\Dp}{\CaD_\kp}
\newcommand{\Dpb}{\CaD_{\kp^b}}
\newcommand{\DV}{\CaD_V}
\newcommand{\DC}{\CaD_\C}
\newcommand{\MF}{\CM_\gl^F}
\newcommand{\IF}{\CI_\gl^F}
\newcommand{\Nl}{\CN_\gl}
\newcommand{\Mlh}{\CM_\gl^\kh}
\newcommand{\Nlh}{\CN_\gl^\kh}
\newcommand{\kgrs}{\kg_{rs}}

\newcommand{\gvtp}{{\gvt_\gP}}

\mainmatter

\section*{Introduction}

Let $G$ be a semi-simple Lie group. An irreducible representation of
$G$ has a character which is an {\sl invariant eigendistribution},
that is a distribution on $G$ which is invariant under the adjoint
action of $G$ and which is an eigenvalue of every biinvariant
differential operator on $G$.  A celebrated theorem of Harish-Chandra
\cite{HC1} says that all invariant eigendistributions are locally
integrable functions on $G$.

After transfer to the Lie algebra $\kg$ of $G$ by the inverse of the
exponential map, an invariant eigendistribution is a solution of a
$\Dg$-module $\CM^F_\gl$ for some $\gl\in\kg^*$. Kashiwara and Hotta
studied in \cite{HOTTA} these $\Dg$-modules $\CM^F_\gl$, in particular
they proved that they are holonomic and, using a modified version of
the result of Harish-Chandra, proved that they are regular
holonomic. This shows in particular that any hyperfunction solution of
a module $\CM^F_\gl$ is a distribution, hence that any invariant
eigenhyperfunction is a distribution.

In \cite{SEKI}, Sekiguchi extended the definition of the modules
$\CM^F_\gl$ to a symmetric pair. A symmetric pair is a decomposition
of a reductive Lie algebra into a direct sum of an even and an odd
part, and the group associated to the even part has an action on the
odd part (see section \ref{sec:def2} for the details). In the diagonal
case where even and odd part are identical, it is the action of
a group on its Lie algebra. Sekiguchi defined a subclass of symmetric
pairs ("nice pairs"), for which he proved a kind of Harish-Chandra
theorem, that is that there is no hyperfunction solution of a module
$\CM^F_\gl$ supported by a hypersurface. He also conjectured that
these modules are regular holonomic.

In \cite{LEVASS} and \cite{LEVASS2}, Levasseur and Stafford give new
proofs of the Harish-Chandra theorem in the original case (the
"diagonal" case) and in the Sekiguchi case ("nice pairs"). In
\cite{GL}, we show that both theorems may be deduced from results on
the roots of the $b$-functions associated to $\CM^F_\gl$.

The aim of this paper is to prove Sekiguchi's conjecture, that is the
regularity of $\CM^F_\gl$, in the general case of symmetric pairs. Our
proof do not use Harish-Chandra's theorem or its generalization, so
we do not need to ask here the pairs to be "nice".

In the first section of the paper we study the regularity of holonomic
$\CaD$-modules. In the definition of Kashiwara-Kawa\"\i \cite{KKREG},
a holonomic $\CaD$-module is regular if it is microlocally regular
along each irreducible component of its characteristic variety. We had
proven in \cite{THESE}, that the microlocal regularity may be
connected to some microcharacteristic variety. We show here that an
analogous result is still true if homogeneity is replaced by some
quasi-homogeneity.

In the second section, we prove Sekiguchi's conjecture in theorem
\ref{thm:conj}. First by standard arguments, we show that outside of
the nilpotent cone, the result may be proved by reduction to a Lie
algebra of lower dimension. Then on the nilpotent cone we use the
results of the first section to show that the module is microlocally
regular along the conormals to the nilpotent orbits.

\section{Bifiltrations of $\CaD$-modules}  \label{chap1}
\subsection{$V$-filtration and microcharacteristic varieties}\label{sec:vfiltr}

In this section, we recall briefly the definitions of the
$V$-filtration and microcharacteristic varieties. Details may be found
in \cite{ENS} (see also \cite{KVAN},\cite{LAUM},\cite{SABB}).

Let $X$ be a complex manifold, $\OX$ be the sheaf of holomorphic
functions on $X$ and $\DX$ be the sheaf of differential operators with
coefficients in $\OX$. Let $Y$ be a submanifold of $X$. The ideal
$\CI_Y$ of holomorphic functions vanishing on $Y$ defines a
filtration of the sheaf $\OX|_Y$ of functions on $X$ defined on a
neighborhood of $Y$ by $F^k_Y\OX=\CI_Y^k$. The associate graduate,
$gr_Y\OX=\bigoplus\CI_Y^k/\CI_Y^{k+1}$ is isomorphic to the sheaf
$\gl_*\CO_{[\tYX]}$ where $\gl:\tYX\to Y$ is the normal bundle to $Y$
in $X$ and $\CO_{[\tYX]}$ the sheaf of holomorphic functions on $\tYX$
which are polynomial in the fibers of $\gl$. For $f$ a function of
$\OX|_Y$ we will denote by $\gs_Y(f)$ its image in $gr_Y\OX$.

If $\CI$ is the ideal of definition of an analytic
subvariety $Z$ of $X$, then $\gs_Y(\CI)=\ensemble{\gs_Y(f)}{f\in \CI}$
is an ideal of $\CO_{[\tYX]}$ which defines the tangent cone to $Z$
along $Y$ \cite{WHIT}.

In local coordinates $(x,t)$ such that $Y=\{t=0\}$, $\CI_Y^k$ is, for 
$k\ge0$, the sheaf of functions
$$f(x,t)=\sum_{\abs{\ga}=k}f_\ga(x,t)t^\ga$$
and if $k$ is maximal with $f\in\CI_Y^k$, we have
$\gs_Y(f)(x,\tt)=\sum_{\abs{\ga}=k}f_\ga(x,0)\tt^\ga$.

Consider now the conormal bundle to $Y$ denoted by $\gL=\TYX$ as a
submanifold of $\TX$. If $f$ is a function on $\TX$, $\gs_\gL(f)$ is a
function  on the normal bundle $T_\gL(T^*X)$. The hamiltonian
isomorphism $TT^*X\simeq T^*T^*X$ associated to the symplectic
structure of $T^*X$ identifies  $T_\gL(T^*X)$ with the the cotangent
bundle $\TL$ and thus  considered $\gs_\gL(f)$ may be considered as a
function on $\TL$. 

\vespa

The sheaf $\DX$ is provided with the filtration by the usual order of
operators denoted by $(\DXi m)_{m\ge0}$ and that we will call the
``usual filtration''. The graduate associated to this filtration is
$gr\DX\simeq \gp_*\CO_{[T^*X]}$ where $\gp:T^*X\to X$ is the cotangent
bundle and $\CO_{[T^*X]}$ is the sheaf of holomorphic functions
polynomial in the fibers of $\gp$. We have also $gr^m\DX\simeq
\gp_*\CO_{[T^*X]}[m]$ where $\CO_{[T^*X]}[m]$ is the sheaf of
holomorphic functions polynomial homogeneous of degree $m$ in the
fibers of $\gp$. If $P$ is a differential operator of $\DX|_Y$, its
principal symbol is a function $\gs(P)$ on $T^*X$ defined in a
neighborhood of $\gL=\TYX$ and $\gs_\gL(\gs(P))$ is a function on
$\TL$ (denoted by $\gs_\gL\{1\}(P)$ in the notations of \cite{ENS}).

The sheaf $\DX|Y$ of differential operators on a neighborhood of $Y$
is also provided with the the $V$-filtration of Kashiwara \cite{KVAN}~:
$$V_k\DX=\{P\in\DX/ \forall j\in\Z,\ \ P\CI_Y^j\subset\CI_Y^{j-k}\}$$
where $\CI_Y^j=\OX$ if $j\le0$.

In local coordinates $(x,t)$, the operators $x_i$ and
$D_{x_i}:=\dsur{x_i}$ have order $0$ for the $V$-filtration while the
operators $t_i$ have order $-1$ and $D_{t_i}:=\dsur{t_i}$ order $+1$.

Remark that the $V$-filtration induces a filtration on $gr\DX\simeq
\gp_*\CO_{[T^*X]}$ which is nothing but the filtration $F_\gL$
associated the conormal bundle $\gL=\TYX$.
In coordinates, $\gL=\ensemble{(x,t,\gx,\gt)\in T^*X}{t=0, \gx=0}$, a
function of $\CO_{[T^*X]}[m]\cap\CI_\gL^{m-k}$ is a function
$f(x,t,\gx,\gt)$ which is polynomial homogeneous of degree $m$ in
$(\gx,\gt)$ and vanishes at order at least $m-k$ on $\{t=0,
\gx=0\}$.

The two filtrations of $\DX$ define a bifiltration
$F_{kj}\DX=\DXi j\cap V_k\DX$. The associated bigraduate is defined by
$gr_F\DX=\oplus gr_F^{kj}\DX$ with
$$gr_F^{kj}\DX=F_{kj}\DX\left/\left(F_{k-1,j}\DX+F_{k,j-1}\DX\right)\right.$$
and is isomorphic to $gr_\gL gr\DX$ that is to the sheaf
$\gp_*\CO_{[T^*\gL]}$ of holomorphic functions on $\TL$
polynomial in the fibers of $\gp:\TL\to Y$. The image of a
differential operator $P$ in this bigraduate will be denoted by
$\gs_\gL\infun(P)$ and may be defined as follows:

If the order of $P$ for the $V$-filtration is equal to the order of
its principal symbol $\gs(P)$ for the induced $V$-filtration then
$\gs_\gL\infun(P)=\gs_\gL(\gs(P))$ and if the order of $\gs(P)$ is
strictly lower then $\gs_\gL\infun(P)=0$.

\vespa 

Let $\CM$ be a coherent $\DX$-module. A good filtration of $\CM$ is a
filtration which is locally finitely generated that is locally of the
form~:
$$\CM_m=\sum_{j=1,\dots,N}\DXi {m+m_j}u_j$$
where $u_1,\dots,u_N$ are (local) sections of $\CM$ and
$m_1,\dots,m_N$ integers.

It is well known that if $(\CM_m)$ is a good filtration of $\CM$,
the associated graduate $gr\CM$ is a coherent $gr\DX$-module and defines the
characteristic variety of $\CM$ which is a subvariety of $T^*X$.
This subvariety is involutive for the canonical symplectic
structure of $T^*X$ and a $\DX$-module is said to be holonomic if
its characteristic variety is lagrangian that is of minimal
dimension. 

In the same way, a good bifiltration of $\CM$ is a bifiltration which
is locally finitely generated. Then the associated bigraduate is a
coherent $gr_F\DX$-module which defines a subvariety
$Ch_\gL\infun(\CM)$ of $T^*\gL$. It is a homogeneous
involutive subvariety of $T^*\gL$ but it is not necessarily lagrangian
even if $\CM$ is holonomic.

If $\CI$ is a coherent ideal of $\DX$ then:
\begin{align*}
\CHM&=\ensemble{\gx\in T^*X}{\forall P\in\CI, \gs(P)(\gx)=0}\\
Ch_\gL\infun(\CM)&=\ensemble{\gz\in \TL}{\forall P\in\CI,
\gs_\gL\infun(P)(\gz)=0}
\end{align*}

\vespa

Regular holonomic $\DX$-modules have been defined by Kashiwara and
Kawa\"\i\ in \cite[Definition 1.1.16.]{KKREG}. A holonomic $\DX$-module
$\CM$ is regular if it has regular singularities along the smooth part
of each irreducible component of its characteristic variety. It is
proved in \cite{KKREG} that the property of regular singularities is
generic, that is it suffices to prove it on a dense open subset of
$\gL$, in particular we may assume that $\gL$ is the conormal bundle
to a smooth subvariety of $X$. The definition of regular singularities
along a smooth lagrangian variety is given in \cite[Definition
1.1.11.]{KKREG} but in this paper, we will use the following
characterization which we proved in \cite[Theorem 3.1.7.]{THESE}:

\begin{prop}\label{charact}
A coherent $\DX$-module has regular singularities along a lagrangian
manifold $\gL$ if and only if $Ch_\gL\infun(\CM)$ is
contained in the zero section of $\TL$.
\end{prop}

\subsection{Weighted $V$-filtration}\label{sec:gbfunct}

The $V$-filtration is associated to the Euler vector field of the
normal bundle $\tYX$ which in coordinates is equal to $\sum
\tt_iD_{\tt_j}$. We want to define a new filtration associated to a
vector field $\sum m_i\tt_iD_{\tt_j}$. As this is not invariant under
coordinate transform, we have first to give an invariant definition.

Let us consider the fiber bundle $p:\tYX\to Y$. The sheaf $\DtYXr$
of relative differential operators is the subsheaf of the sheaf
$\DtYX$ of differential operators on $\tYX$ commuting with all
functions of $p^{-1}\OY$. A differential operator $P$ on $\tYX$
is homogeneous of degree $0$ if for any function $f$ homogeneous
of degree $k$ in the fibers of $p$, $Pf$ is homogeneous of degree
$k$.

In particular, a vector field $\tgh$ on $\tYX$ which is a relative
differential operator homogeneous of degree $0$ defines a morphism
from the set of homogeneous functions of degree $1$ into itself which
commutes with the action of $p^{-1}\OY$, that is a section of
$$\HOM_{p^{-1}\OY}\left(\CO_\tYX[1],\CO_\tYX[1]\right).$$
Let $(x,t)$ be coordinates of $X$ such that $Y=\ensemble{(x,t)\in
X}{t=0}$. Let $(x,\tt)$ be the corresponding coordinates of
$\tYX$. Then $\tgh$ is written as~:
$$\tgh=\sum a_{ij}(x)\tt_iD_{\tt_j}$$ and the matrix $A=(a_{ij}(x))$
is the matrix of the associated endomorphism of $\CO_\tYX[1]$ which is
a locally free $p^{-1}\OY$-module of rank $d=codim_XY$. Its
conjugation class is thus independent of the choice of coordinates
$(x,t)$. When the morphism is the identity, $\tgh$ is by definition
the Euler vector field of $\tYX$.

\begin{defi}
A vector field $\tgh$ on $\tYX$ is {\sl definite positive} if it is a
relative differential operator homogeneous of degree $0$ whose
eigenvalues are strictly positive rational numbers and which is
locally diagonalizable as an endomorphism of $\CO_\tYX[1]$.
\end{defi}

A structure of \textsl{local fiber bundle of $X$ over $Y$} is  an
analytic isomorphism between a neighborhood of $Y$ in $X$ and a
neighborhood of $Y$ in $\tYX$. For example a local system of
coordinates defines such an isomorphism. 

\begin{defi}
A vector field $\gh$ on $X$ is \textsl{definite positive with respect
to $Y$} if:

(i)\ $\gh$ is of degree $0$ for the $V$-filtration associated to $Y$
and the image $\gs_Y(\gh)$ of $\gh$ in $gr^0_V\DX$ is {\sl definite
positive} as a vector field on $\tYX$.

(ii)\ There is a structure of local fiber bundle of $X$ over $Y$ which
identifies $\gh$ and $\gs_Y(\gh)$.
\end{defi}

It is proved in \cite[proposition 5.2.2]{ENS} that if $\gs_Y(\gh)$ is
the Euler vector field of $T_YX$ the condition (ii) is always
satisfied and the local fiber bundle structure of $X$ over $Y$ is
unique for a given $\gh$, but this is not true in general.

\vespa

We will now assume that $X$ is provided with such a vector field
$\gh$. Let $\gb=a/b$ the rational number with minimum positive
integers $a$ and $b$ such that the eigenvalues of $\gb^{-1}\gh$ are
positive relatively prime integers. Let $\DX[k]$ be the sheaf of
differential operators $Q$ satisfying the equation $[Q,\gh]=\gb kQ$
and let $V^\gh_k\DX$ be the sheaf of differential operators $Q$ which
are equal to a series $Q=\sum_{l\le k}Q_l$ with $Q_l$ in $\DX[l]$ for
each $l\in\Z$.

By definition of a definite positive vector field, we may find local
coordinates $(x,t)$ such that $\gh=\sum m_it_iD_{t_i}$ and we may
assume that the $m_i$ are relatively prime integers after
multiplication of $\gh$ by $\gb^{-1}$. In this situation, the operators
$x_j$ and $D_{x_j}$ have order $0$ while the operators $t_i$ have
order $-m_i$ and $D_{t_i}$ order $+m_i$. This shows in particular that
any monomial $x^\ga t^\gb D_x^\gg D_t^\gd$ is in some $\DX[k]$ and
thus that $\DX$ is the union of all $V^\gh_k\DX$. This defines a
filtration $V^\gh$ of the sheaf of rings $\DX$.

The principal symbol of $[Q,\gh]$ is the Poisson bracket
$\{\gs(P),\gs(\gh)\}$ which is equal to $H_\gh(\gs(P))$ where  $H_\gh$
is a vector field on $T^*X$, the Hamiltonian of $\gh$.
The $V^\gh$-filtration on $\DX$ induces a filtration on the graduate
of $\DX$ that is on $\CO_{[T^*X]}$.  A function $f$ of $\CO_{[T^*X]}$ will be in
$V^\gh_k\CO_{[T^*X]}$ if it is a series of functions $f_l$ for $l\ge
k$ with $H_\gh f = -lf$. In this case we set $\gs^\gh_k(f)=f_k$.

\vespa

We are now in a situation analog to that of section \ref{sec:vfiltr}
with two filtrations on $\DX$, the usual filtration and the
$V^\gh$-filtration. The sheaf $\DX$ is thus provided with a
bifiltration by $F^\gh_{kj}\DX=\DXi j\cap V^\gh_k\DX$ and this defines a
symbol $\gs^\gh\infun(P)$ which is a function on $\TX$. By definition,
$\gs^\gh\infun(P)$ is equal to $\gs^\gh_k(\gs(P))$ where $k$ is the
order of $P$ for the $V^\gh$-filtration. This symbol is thus equal to
$0$ if the order of $\gs(P)$ is strictly less than $k$.

If $\CM$ is a
coherent $\DX$-module, we define a good bifiltration and a
microcharacteristic variety $Ch^\gh\infun(\CM)$. If $\CM=\DX/\CI$ we
will have:
$$Ch^\gh\infun(\CM) =\ensemble{\gz\in\TX}{\forall
P\in\CI,\quad \gs^\gh\infun(P)(\gz)=0}$$

The difference with the previous situation is the local identification
of $\tYX$ with $X$ which defines isomorphisms $T^*\TYX\simeq
T^*\tYX\simeq T^*X$ and make $\gs^\gh(\gs(P))$ a function on
$T^*X$. Especially, if $\tgh$ is the Euler vector field of $\tYX$ and
$\gh$ a vector field on $X$ with $\gs_V(\gh)=\tgh$, the definitions of
this section coincide with the definitions of the previous one except
for this identification.

\subsection{Direct image of $V$-filtration}

Let $\gf:Y\to X$ be a morphism of complex analytic manifolds. A vector
field $u$ on $Y$ is said to be \textsl{tangent to the fibers of $\gf$}
if $u(f\circ\gf)=0$ for all $f$ in $\OX$. A differential operator $P$
is said to be \textsl{invariant under $\gf$} if there exists a
$\C$-endomorphism $A$ of $\OX$ such that $P(f\circ\gf)=A(f)\circ\gf$
for all $f$ in $\OX$. If we assume from now that $\gf$ has a dense
range in $X$, $A$ is uniquely determined by $P$ and is a differential
operator on $X$. We will denote by $A=\gf_*(P)$ the image of $P$ in
$\DX$ under this ring homomorphism.

Let $Z$ be a submanifold of $Y$ and $T$ a submanifold
of $X$. Let $\gh$ be a vector field on $Y$ invariant under $\gf$. We
assume that $\gh$ is definite positive with respect to $Z$ and that
$\gh'=\gf_*(\gh)$ is definite positive with respect to $T$. We also
multiply $\gh$ by an integer so that its eigenvalues and those of
$\gh'$ are integers.

\begin{exem}\label{ex:}
Let $Y$ be a complex vector space and $\gf:Y\to X=\C^d$ given by
$\gf=(\gf_1,\dots,\gf_d)$ where $\gf_1,\dots,\gf_d$ are holomorphic
functions on $Y$ homogeneous of degree $m_1,\dots,m_d$. Let $Z=\{0\}$ and
$\gh$ be the Euler vector field of $Y$, so that the $V^\gh$-filtration is
the $V$-filtration along $\{0\}$. Then $\gh'=\gf_*(\gh)$ is equal to
$\sum m_it_iD_{t_i}$ on $X$ and is definite positive with respect to
$\{0\}$. Remark that we do not assume that $\gf$ is defined in a
neighborhood of $Z$.
\end{exem}

In the general case, we can choose local coordinates $(y,t)$ on $X$ so
that $\gh'=\sum m_jt_jD_{t_j}$, then the map $\gf$ is given by
$y_i=\gf_i(x)$ and $t_j=\gy_j(x)$ where the functions $\gf_i(x)$ is
homogeneous of degee $0$ for $\gh$ while the function $\gy_j(x)$ is
homogeneous of degee $m_j$ for $\gh$.

\vespa

The sheaf $\DYX=\OY\ox_{\gf^{-1}\OX}\gf^{-1}\DX$ is a
($\DY$,$\gf^{-1}\DX$)-bimodule with a canonical section $1\ox1$
denoted by $\uYX$.  If we choose coordinates $(x_1,\dots,x_n)$ of $X$
and coordinates $(y_1,\dots,y_p)$ of $Y$ and if
$\gf=(\gf_1,\dots,\gf_p)$, then the sections of $\DYX$ are represented
by finite sums $\sum f_\ga(y)\ox D_x^\ga$ and the left action of $\DY$
is given by
\[D_{y_i}\left(\sum_\ga f_\ga(y)\ox D_x^\ga\right) =
\sum_{\ga}\frac{\partial f_\ga}{\partial y_i}(y)\ox D_x^\ga
+\sum_{\ga,j}f_\ga(y)\frac{\partial\gf_j}{\partial y_i}(y)\ox D_{x_j}D_x^\ga\]

If $\CN$ is a coherent $\DX$-module, its inverse image under
$\gf$ is the $\DY$-module
$\gf^*\CN=\DYX\ox_{\gf^{-1}\DX}\gf^{-1}\CN$. In general,
$\gf^*\CN$ is not coherent but if $\CN$ is holonomic, $\gf^*\CN$
is holonomic (hence coherent).

Let $\DYX[k]$ be the set of sections satisfying
$\gh.u-u.\gf_*\gh=-\gb\gb'ku$ where $\gb$ (resp. $\gb'$) is the g.c.d. of
the eigenvalues of $\gh$ (resp. $\gf_*\gh$). (We may
assume that $\gb=1$ or $\gb'=1$ but not both in general). We define
$V_k\DYX$ as the subsheaf of  $\DYX$ of the sections which may be
written as series $\sum_{l\ge k}u_l$ with $u_l$ in
$\DYX[l]$. Remark that $\uYX$ satisfies $\gh.\uYX=\uYX.\gf_*\gh$
hence is of order $0$.

If  $\CN$ is a coherent $\DX$-module provided with a 
$V^{\gh'}$-filtration we define a filtration on its inverse image by:

$$V^\gh_k\gf^*\CN=\sum_{k=\gb' i+\gb j}V_i\DYX\ox_{\gf^{-1}V^{\gh'}_0\DX}
\gf^{-1}V^{\gh'}_j\CN$$

The sheaf $\DYX$ is also provided with a filtration $(\DYX)_j$ induced
by the usual filtration of $\DX$ hence of a bifiltration
$F^\gh\DYX$. If $\CN$ is bi-filtrated, we define in the same way a
bifiltration on $\gf_*\CN$.

\begin{prop}\label{prop:transfer}
Let $\CI$ be an ideal of $\DY$ which is generated by all the
vector fields tangent to the fibers of $\gf$ and by a finite set
$(P_1,\dots,P_l)$ of differential operators invariant under
$\gf$. Let $\CJ$ be the ideal of $\DX$ generated by
$(\gf_*(P_1),\dots,\gf_*(P_l))$. Let $\CM=\DY/\CI$ and
$\CN=\DX/\CJ$ and put on $\CM$ and $\CN$ the bifiltrations induced by
$F^\gh\DY$ and $F^{\gh'}\DX$

Then, there exists a canonical morphism of $\DY$-modules
$\CM\to\gf^*\CN$ which is  a morphism of bi-filtrated
$F^\gh\DY$-modules and an isomorphism at the points where
$\gf$ is a submersion.
\end{prop}

\begin{proof}
There is a canonical morphism $\DY\to\DYX$ given by
$P\mapsto P.\uYX$. The vector fields tangent to the fibers cancel
$\DYX$ and a differential operator invariant under $\gf$ satisfy
$P.\uYX=\uYX.\gf_*(P)$ hence this morphism defines a morphism
$\CM\to\gf_*\CN$ which is a morphism of left $V^\gh\DY$-modules by the
definitions.

In a neighborhood of a point where $\gf$ is a submersion, we may
choose local coordinates $(x_1,\dots,x_{p},y_1,\dots,y_{n-p})$ such
that $\gf(x,y)=x$. Then $\DYX$ is the sheaf of operators $P(x,y,D_x)$,
the vector fields tangent to the fibers are generated by
$D_{y_1},\dots,D_{y_{n-p}}$ and the differential operators invariant
under $\gf$ are of the form $P(x,D_x)$ modulo $(D_{y_i})$, so
$\CM\to\gf_*\CN$ is an isomorphism.
\end{proof}

Let $S=\gf^{-1}(T)$ and $x$ be a point of $S$ where $\gf$ is a
submersion. In a neighborhood of $x$, $Y$ is isomorphic to $X\x S$ and
if we fix such an isomorphism, $\gh'$ which is a vector field on $X$
may be considered as a vector field on $Y$, definite positive
relatively to $S$. Remark that $\gh'$ differ from $\gh$ by a vector
field tangent to $\gf$. Then proposition \ref{prop:transfer} gives:

\begin{coro}\label{cor:transfer}
The microcharacteristic variety $Ch^\gh\infun(\CM)$ is equal to
$Ch^{\gh'}\infun(\CM)$ in a neighborhood of $x$.
\end{coro}

\subsection{Weighted $V$-filtration and regularity}

\begin{defi}\label{def:wrs}
Let $Z$ be a submanifold of $X$ and $\gh$ be a vector field which is
definite positive with respect to $Z$.  A coherent $\DX$-module has
$\gh${\sl -weighted regular singularities} along the lagrangian
manifold $\gL=T^*_ZX$ if there is a dense open subset $\gO$ of $\gL$ such
that $Ch^\gh\infun(\CM)\subset\gL$ in a neighborhood of $\gO$.
\end{defi}

If $\gs_Z(\gh)$ is the Euler vector field of $T_ZX$, proposition
\ref{charact} shows that this definition coincide with the definition
of Kashiwara-Kawa\"\i.

\vespa

Let $X=\C^n$ with coordinates $(x_1,\dots,x_{n-p},t_1,\dots,t_p)$ and
$Z=\{t=0\}$, let $Y=\C^n$ with coordinates
$(x_1,\dots,x_{n-p},y_1,\dots,y_p)$ and $Z'=\{y=0\}$. Let
$m_1,\dots,m_p$ be strictly positive integers , we define the map
$\gf:Y\to X$ by $\gf(x,y)=(x,y_1^{m_1},\dots,y_p^{m_p})$ and the
vector field $\gh=\sum_{i=1\dots p}m_it_iD_{t_i}$.

\begin{lemm} \label{lem:wr}
Let $\CM$ be a holonomic $\DX$-module with $\gh$-weighted
regular singularities along $T^*_ZX$, then $\gf^*\CM$ is a holonomic
$\DY$-module with regular singularities along $T^*_{Z'}Y$.
\end{lemm}

\begin{proof}
We may assume that $\CM$ is equal to $\DX/\CI$ for some coherent ideal
$\CI$ of $\CM$. The inverse image of $\CM$ by $\gf$ is, by definition:
$$\gf^*\CM=\DYX\otimes_{\gf^{-1}\DX}\gf^{-1}\CM=\DYX/\DYX\CI$$

The sections of $\DYX$ are represented by $P(x,y,D_x,D_t)=\sum
a_{\ga\gb}(x,y)D_x^\ga D_t^\gb$ and we define the filtration
$V^\gh\DYX$ in the same way than in the previous section. For this
filtration $x^\gg y^\gd D_x^\ga D_t^\gb$ is of order $\scal
m\gb\!\!-|\gd|$. We also define the usual filtration on $\DYX$, that is the
filtration by the order in $(D_x,D_t)$. In this way, $\DYX$ is
provided with a bifiltration $F^\gh\DYX$ which is compatible with the
bifiltration $F^\gh\DX$, that is an operator $P$ of  $F_{kl}^\gh\DX$
sends  $F_{ij}^\gh\DYX$ into $F_{i+k,j+l}^\gh\DYX$.

Let $\DYX[N]$ be the sub-$\DY$-module of $\DYX$ generated by $D_t^\gb$
for $|\gb|\le N$. If $\CM$ is holonomic, $\gf^*\CM$ is holonomic hence
coherent. The images of the morphisms $\DYX[N]\to\gf^*\CM$ make an
increasing sequence of coherent submodules of $\gf^*\CM$ which is
therefore stationary, so there exists some $N_0$ such that
$\DYX[N]\to\gf^*\CM$ is surjective for all $N\ge N_0$. The
bifiltration induced by $F^\gh\DYX$ on $\DYX[N]$ is a good
$F\DY$-filtration which induces a good filtration on $\gf^*\CM$ if
$N\ge N_0$, we will denote it by $F[N]\gf^*\CM$.

The associate graduate is denoted by $gr[N]\gf^*\CM$ and, as $F[N]$ is
a good bifiltration, the analytic cycle of $T^*Y$ associated to
$gr[N]\gf^*\CM$ is independent of $N$ \cite[Prop 3.2.3.]{ENS}. For
$N\ge N_0$, the canonical morphism $gr[N_0]\gf^*\CM \to gr[N]\gf^*\CM$
induces an isomorphism on the associated cycles hence
$gr[N_0]\gf^*\CM$ and $gr[N]\gf^*\CM$ have the same support and the
kernel and cokernel of the morphism have a support of dimension
strictly lower.

An operator $P$ of $F_{kl}^\gh\DX$ sends $F^\gh_{ij}\DYX[N_0]$ into
$F^\gh_{i+k,j+l}\DYX[N_0+l]$. If $P$ annihilates a section $u$ of
$F_{ij}[N_0]\gf^*\CM$, its class in $gr_{kl}\DX$ that is the function
$\gs^\gh\infun(P)$ annihilates the image of $u$ in $gr[N+l]\gf^*\CM$.
Let $\gz$ be a point of $\gL= T^*_ZX$ such that
$Ch^\gh\infun(\CM)\subset T^*_ZX$ in a neighborhood of $\gz$. By the
hypothesis, there is a dense open subset $\gO$ of such points in
$\gL$. There is a differential operator $P$ which annihilates $u$ and
such that $\gs^\gh\infun(P)=t_1^M\gm$ where $\gm$ is a function
invertible at $\gz$. Hence there exists some $l$ such that the image
of $u$ in $gr[N+l]\gf^*\CM$ is annihilated by $t_1^M=y_1^{Mm_1}$ hence
is supported by $y_1=0$. As $gr[N_0]\gf^*\CM$ is finitely generated,
there exists some $N_1\ge N_0$ such that the image of $gr[N_0]\gf^*\CM$
in $gr[N_1]\gf^*\CM$ is contained in $y_1=0$.

We can do the same for the other equations of $T^*_{Z'}Y$ and show
that there exists some $N_2\ge N_0$ such that the image of
$gr[N_0]\gf^*\CM$ in $gr[N_2]\gf^*\CM$ is contained in $T^*_{Z'}Y$.
This shows that $gr[N_0]\gf^*\CM$ is supported by the union of
$T^*_{Z'}Y$ and of a set $W$ of dimension strictly lower than the
dimension of $T^*_{Z'}Y$. But we know that this support is involutive
hence all its component have a dimension at least that dimension, so
$gr[N_0]\gf^*\CM$ is supported in $T^*_{Z'}Y$ in a neighborhood of
$\gf^{-1}(\gz)$. By definition $gr[N_0]\gf^*\CM$ is equal to
$Ch_{T^*_{Z'}Y}\infun(\gf^*\CM)$, hence $\gf^*\CM$ has regular singularities
along $T^*_{Z'}Y$.
\end{proof}

\begin{theo}\label{thm:wregular}
Let $X$ be a complex manifold, $\gp:\TX\to X$ the projection,
$Z$ a submanifold of $X$ and $\gh$ a vector field on $X$ which is
definite positive with respect to $Z$. Let $\CM$ be a holonomic
$\DX$-module. We assume that:
\begin{enumerate}
\item$\CM$ is a regular holonomic $\DX$-module on $X-Z$,\label{enum1}
\item$\CM$ has $\gh$-weighted regular singularities along $T^*_ZX$,\label{enum2}
\item The dimension of $\CHM\cap T^*_ZX$ is equal to the dimension of $X$. \label{enum3}
\end{enumerate}
Then $\CM$ is a regular holonomic $\DX$-module.
\end{theo}

\begin{proof} We fix local coordinates $(x_1,\dots,x_{n-p},t_1,\dots,t_p)$ of $X$ so
that $Z=\{t=0\}$ and $\gh=\sum_{i=1\dots p}m_it_iD_{t_i}$. We define a map $\gf:Y\to
X$ by $\gf(x,y)=(x,y_1^{m_1},\dots,y_p^{m_p})$ where $Y$ is a
neighborhood of $0$ in $\C^n$. If $Z'$ is the set $\{y=0\}$, lemma
\ref{lem:wr} shows that $\gf^*\CM$ has regular singularities along
$T^*_{Z'}Y$.

The third condition means that the characteristic
variety of $\CM$ has no irreducible component contained in
$\gp^{-1}(Z)$ except $T^*_ZX$. The same is true for $\gf^*\CM$ on
$Z'$. This may be proved as in lemma \ref{lem:wr}
but with the usual filtration replacing the
bifiltration. This may also be proved easily with the definition of
the characteristic variety in terms of microdifferential operators. 

By hypothesis, $\CM$ is regular on $X\setminus Z$ hence by \cite[Cor
5.4.8.]{KKREG} $\gf^*\CM$ is regular holonomic on $Y\setminus Z'$. So,
$\gf^*\CM$ has regular singularities along each irreducible component
of its characteristic variety, hence by definition, it is a regular
holonomic $\DY$-module.

Then by \cite[theorem 6.2.1.]{KKREG}, the direct image $\gf_*\gf^*\CM$
is a regular holonomic $\DX$-module. By definition
$$\gf_*\gf^*\CM
=\R\gf_*(\DXY\otimes_\DY^\L\DYX\otimes_{\gf^{-1}\DY}^\L\gf^{-1}\CM)$$
and the morphism $\DX\to\DXY\otimes_\DY\DYX$ is injective hence$\CM$
is a submodule of $\gf_*\gf^*\CM$ hence a regular holonomic $\DX$-module.
\end{proof}

The following corollary is the generalization of the definition of
regular holonomic $\CaD$-modules and of proposition \ref{charact}. It
is proved from the previous theorem by descending induction on the
dimension of the strata.

\begin{coro} Let $\CM$ be a holonomic $\DX$-module. Assume that there
is a stratification $X=\bigcup X_\ga$ such that $\CHM\subset\bigcup
T^*_{X_\ga}X$ and for each $\ga$ there is a vector field $\gh_\ga$
positive definite along $X_\ga$ such that $\CM$ has $\gh_\ga$-weighted
regular singularities along $T^*_{X_\ga}X$.

Then $\CM$ is a regular holonomic $\DX$-module.
\end{coro}

\section{Symmetric pairs}  \label{chap2}
\subsection{Definitions}\label{sec:def2}

Let us briefly recall what is a symmetric pair. For the details we
refer to \cite{SEKI} and \cite{LEVASS2}. Let $G$ be a connected
complex reductive algebraic group with Lie algebra $\kg$. Fix a
non-degenerate, $G$-invariant symmetric bilinear form $\gk$ on the
reductive Lie algebra $\kg$ such that $\gk$ is the Killing form on the
semi-simple Lie algebra $[\kg,\kg]$. Fix an involutive automorphism
$\gvt$ of $\kg$ preserving $\gk$ and set $\kk=Ker(\gvt-I)$,
$\kp=Ker(\gvt+I)$. Then $\kg=\kk\oplus\kp$ and the pair $(\kg,\kk)$ or
$(\kg,\gvt)$ is called a symmetric pair. Recall that $\kk$ and $\kp$
are orthogonal with respect to $\gk$ and that $\kk$ is a reductive Lie
subalgebra of $\kg$. Denote by $K$ the connected reductive subgroup of
$G$ with Lie algebra $\kk$. The group $K$ acts on $\kp$ via the
adjoint action.

Let $\kp^*$ be the dual of $\kp$, $\CO(\kp)=S(\kp^*)$ the
ring of regular functions on $\kp$ ($S(\kp^*)$ is the symmetric
algebra), $\CO(\kp^*)=S(\kp)$ the ring of regular functions on $\kp^*$
and $\CaD(\kp)$ the ring of differential operators on $\kp$ with
coefficients in $\CO(\kp)$. The ring of functions $\CO(\kp)$ is
naturally embedded in $\CaD(\kp)$ and we embed $\CO(\kp^*)=S(\kp)$ in
$\CaD(\kp)$ as differential operators with constant coefficients. That
is we associate to an element $u$ of the vector space $\kg$ the
derivation in the direction of $u$
$$D_u(f)(x)=\frac{d}{dt}f(x+tu)|_{t=0}$$
and we extend to the symmetric algebra $S(\kp)$. Remark that this
embedding is compatible with the filtration by the degree in $S(\kp)$
and the filtration by the order in $\CaD(\kp)$.

Notice that $K$ has an induced action on $S(\kp)$, $S(\kp^*)$ and
$\CaD(\kp)$ and we have natural embeddings of the invariant subrings
$S(\kp)^K\subset\CaD(\kp)^K$ and $S(\kp^*)^K\subset\CaD(\kp)^K$. The
ring $S(\kp)^K$ is equal to the ring of polynomials
$\C[p_1,\dots,p_r]$ for some $p_1,\dots,p_r$ in  $S(\kp)^K$ and in the
same way  $S(\kp^*)^K$ is equal to a ring of polynomials
$\C[q_1,\dots,q_r]$ \cite{KOSR}.

The differential of the action of  $K$ on $\kp$ induces a Lie
algebra homomorphism $\gt:\kk\to \mathrm{Der}\ S(\kp^*)$ hence an
embedding $\gt:\kk\to\CaD(\kp)$ defined by
$$(\gt(a).f)(v)=\frac{d}{dt}f(e^{-ta}.v)|_{t=0},\quad 
\mathrm{for}\ a\in\kk,\ f\in \CO(\kp),\ v\in\kp$$

As a section of the tangent bundle, $\gt(A)$ is the map $\kp\to
T\kp=\kp\x\kp$ given by $\gt(A)(X)=(X,[X,A])$.

We denote by $\bN(\kp)$ the nilpotent cone of $\kp$, that is the set
of nilpotent elements of $\kg$ which lie in $\kp$, it is also the
subvariety of $\kp$ defined by the set of $K$-invariant functions
$S(\kp^*)^K$. In the same way we consider the nilpotent cone
$\bN(\kp^*)$ which is the subvariety of $\kp^*$ defined by $S(\kp)^K$.
An important result is that the nilpotent cone $\bN(\kp)$ is a finite
union of $K$-orbits \cite[theorem 2]{KOSR}.

The cotangent bundle $T^*\kp$ is equal to $\kp\x\kp^*$. The
non-degenerate form $\gk$ on $\kg$ defines a non-degenerate symmetric
bilinear form on $\kp$ and an isomorphism $\kp\simeq\kp^*$. We
identify $T^*\kp=\kp\x\kp^*\simeq\kp\x\kp$. Let
$\CC(\kp)=\ensemble{(x,y)\in\kp\x\kp}{[x,y]=0}$, then the dimension of
$(\kp\x\bN(\kp))\cap\CC(\kp)$ is equal to the dimension of $\kp$
\cite[lemma 2.2.]{LEVASS2}.

The characteristic variety of $\Dp/\Dp\gt(\kk)$ is equal to
$\CC(\kp)$. Let $F$ be an ideal of finite codimension of $S(\kp)^K$,
its graduate is a power of $S(\kp)^K$ hence the characteristic variety
of the $\Dp$-module $\Dp/\Dp F$ is $\kp\x\bN(\kp)$. Finally, if $\CI$
be the left ideal of $\Dp$ generated by $F$ and $\gt(\kk)$, the
characteristic variety of $\CM_F=\Dp/\CI$ is contained in
$(\kp\x\bN(\kp))\cap\CC(\kp)$ hence $\CM_F$ is a holonomic
$\Dp$-module.

As a special case, we have the diagonal case where $G=G_1\x G_1$ with
$\gvt(x,y)=(y,x)$ for some reductive group $G_1$. Thus
$(\kg,\kk)=(\kg_1\oplus\kg_1,\kg_1)$ and $K=G_1$ with its adjoint
action on $\kp=\kg_1$. Let $\gl\in\kp^*$ and
$F_\gl=\ensemble{P-P(\gl)}{P\in S(\kp)^K}$, then the corresponding
module $\CM^F_\gl=\Dp/\Dp\gt(\kk)+\Dp F_\gl$ is the module of
Kashiwara-Hotta \cite{HOTTA}.

\subsection{The conjecture of Sekiguchi}

\begin{theo}\label{thm:conj}
Let $F$ be an ideal of finite codimension of $S(\kp)^K$ and
$\CM_F=\Dp/\CI$ where $\CI$ is the left ideal of $\Dp$ generated by $F$
and $\gt(\kk)$.

Then $\CM_F$ is a \textsl{regular} holonomic $\Dp$-module.
\end{theo}

The proof of this theorem will be made in several steps. First we will
reduce to the semi-simple case (lemma \ref{lem:ss}), then prove
by induction on the dimension of the Lie algebra, that the result is
true outside of the nilpotent cone (lemma \ref{lem:closed}) and the
key point of the proof is the case of a nilpotent orbit (lemma
\ref{lem:nilorbit}).

\begin{lemm}\label{lem:fin}
Let $Y$ be a complex manifold and $X=Y\x\C$. Let $P(t,D_t)$ be a
differential operator on $\C$ with principal symbol independent of $t$
and $\CI$ be a coherent ideal of $\DX$ which contains $P$.

Let $\CM_Y$ be the inverse image of $\CM=\DX/\CI$ on $Y$ by the
immersion $Y\to X$, then $\CM$ is isomorphic to the inverse image of
$\CM_Y$ by the projection $q:X\to Y$, that is
$$\CM=\DXtY\otimes_{q^{-1}\DY}q^{-1}\CM_Y=\CM_Y\widehat\ox\CO_\C$$
In particular, $\CM$ is regular holonomic if and only if $\CM_Y$ is
regular holonomic.
\end{lemm}

\begin{proof}
This lemma is a (very) special case of \cite[theorem 5.3.1. ch
II]{SKK}. The first step is to prove that $\DC/\DC P$ is isomorphic to
$(\DC/\DC D_t)^N$. The proof is the same than that of \cite[theorem
5.2.1. ch II]{SKK}, but as there is only one variable, the proof is
very simple and use only functions instead of differential operators
of infinite order. Then we can follow the proof of \cite{SKK} but with
finite order operators instead of infinite order operators.

Remark that if $P$ were a differential operator in several variables,
for example, $P=D_t^2+D_x$, this result would be true only with the
sheaf $\mathcal D_X^\infty$ of differential operators with infinite order.

As $X=Y\x\C$, the inverse image of $\CM_Y$ by $q$
is isomorphic to the external product of $\CaD$-modules
$\CM_Y\widehat\ox\CO_\C$.
\end{proof}

Assume that $\kp=\kp_0\oplus\kp_1$, the action of $K$ on $\kp_0$ being
trivial. Then $S(\kp)^K=S(\kp_0)\ox S(\kp_1)^K$, this defines a
morphism $\gd:S(\kp)^K\to S(\kp_1)^K$ by restriction and $F_1=\gd(F)$
is an ideal of finite codimension of $S(\kp_1)^K$. Let
$\CM_{F_1}=\CaD_{\kp_1}/\CI_1$ where $\CI_1$ is the ideal of
$\CaD_{\kp_1}$ generated by $\gt_{\kp_1}(\kk)$ and $F_1$.

\begin{lemm}\label{lem:ss}
(1) The module $\CM_F$ is isomorphic to
$\CO_{\kp_0}\widehat\ox(\CM_F)_{\kp_1}$ where $(\CM_F)_{\kp_1}$ is the
restriction of $\CM_F$ to $\kp_1$.

(2) $(\CM_F)_{\kp_1}$ (hence $\CM_F$) is regular if
$\CM_{F_1}$ is regular.
\end{lemm}

\begin{proof}
By induction on the dimension of $\kp_0$, we may assume that $\kp_0=\C$
and choose linear coordinates $(x,t)$ of $\kp$ such that
$\kp_0=\ensemble{(x,t)\in\kp}{x=0}$.  The action of $K$ is trivial on
$\kp_0$ hence $S(\kp)^K$ contains $S(\kp_0)$ and as $F$ is finite
codimensional in $S(\kp)^K$ it contains a polynomial in
$D_{t}$. Lemma \ref{lem:fin} shows the first part of the
lemma.

We assume now that $\CM_{F_1}$ is regular.  Recall that
$(\CM_F)_{\kp_1}=\CM_F/t\CM_F$ is a holonomic $\CaD_{\kp_1}$-module generated by
the classes of $1,\dots,D_t^{m-1}$. Let $\CM'$ be the submodule of
$(\CM_F)_{\kp_1}$ generated by the class $u$ of $D_t^{m-1}$. The
vector fields of $\gt(\kk)$ are independent of $(t,D_t)$ hence $u$ is
annihilated by $\gt(\kk)$. If $P$ is an element of $F$, as an operator
of $\Dp$ it is equal to $\gd(P)+AD_t$ hence $\gd(P)$ annihilates
$u$. So $u$ is annihilated by $\gt(\kk)$ and by $F_1$ and $\CM'$ is a
quotient of $\CM_{F_1}$. So $\CM'$ is regular.

Consider now $\CM''$ which is the submodule of $\CM$ generated by the
classes $D_t^{m-1}$ and $D_t^{m-2}$. The quotient
$\CM''/\CM'$ is generated by the class $v$ of $D_t^{m-2}$ which
is annihilated by $\gt(\kk)$ and by $F_1$, so it is regular. We have an
exact sequence $\exacte{\CM'}{\CM''}{\CM''/\CM'}$ where two terms
are regular hence $\CM''$ is regular. Continuing the same argument, we
get that  $(\CM_F)_{\kp_1}$ is regular.
\end{proof}

Let $b$ be a semisimple element of $\kp$. Then
$\kp=\kp^b\oplus[\kk,b]$ and $\kg^b=\kk^b\oplus\kp^b$ defines a
symmetric pair. Let $\gd$ be the restriction map $\gd:S(\kp)^K\to
S(\kp^b)^{K^b}$, this map is injective and if $F$ is an ideal of
finite codimension of $S(\kp)^K$ then $\gd(F)$ is an ideal of finite
codimension of $S(\kp^b)^{K^b}$ \cite[lemma 19]{HC2}. Let $\CI_b$ be
the left ideal of $\Dpb$ generated by $\gd(F)$ and $\gt(\kk^b)$ and
$\CM_b=\Dpb/\CI_b$.

\begin{lemm}\label{lem:closed}
In a neighborhood of $b$, $\CM_F$ is isomorphic to the external
product of the holomorphic functions on the orbit of $b$ by a quotient
of $\CM_b$. In particular, $\CM_F$ is regular if $\CM_b$ is
regular.
\end{lemm}

\begin{proof}
Let $V$ be a linear subspace of $\kk$ such that
$\kk=V\oplus\kk^b$. The map $f:V\x\kp^b\to\kp$ given by
$f(y,Z)=exp(y).Z$ is a local isomorphism. If $(x_1,\dots,x_{n-r})$ are
linear coordinates of $V$ and $(t_1,\dots,t_r)$ are linear coordinates
of $\kp^b$, the map $f$ defines local coordinates
$(x_1,\dots,x_{n-r},t_1,\dots,t_r)$ of $\kp$ in a neighborhood of $b$.
Lemma 3.7 of  \cite{SEKI} shows that
in these coordinates, the orbit $Kb$ is $\ensemble{(x,t)}{t=0}$,
$\kp^b=\ensemble{(x,t)}{x=0}$ and the differential operators
$D_{x_1},\dots,D_{x_{n-r}}$ belong to $\gt(\kk)$. Hence $\CM$ is the
product of $\CO_{Kb}$ by a module $\CN$. 

If $Z$ is an element of $\kk^b$, $\gt_\kp(Z)$ is by definition the
vector field on $\kp$ with value $[Z,A]$ at a point $A$ of $\kp$. The
value of $\gt_{\kp^b}(Z)$ at a point $A$ of $\kp^b$ is the projection
of $[Z,A]$ on $\kp^b$, hence $\gt_{\kp^b}(\kk^b)$ is equal to
$\gt_\kp(\kk)$ modulo $D_{x_1},\dots,D_{x_{n-r}}$. On the other hand,
let $P\in F$, as the coordinates $(t_1,\dots,t_r)$ are linear
coordinates of $\kp^b$, the value of $P$ on a function of $t$ is the
restriction of $P$ to $S(\kp^b)^{K^b}$. Hence $\CN$ is a quotient of
$\CM_b$.
\end{proof}

\begin{lemm}\label{lem:pp}
Let $\gL$ be the conormal to $0$ in $\kp$. The microcharacteristic
variety $Ch_\gL\infun(\CM_F)$ is contained in
$(\kp\x\bN(\kp))\cap\CC(\kp)$.
\end{lemm}

\begin{proof}
Let $E$ be the Euler vector field of the vector space $\kp$. It is
clear on the definition, that the vector fields of $\gt(\kk)$ preserve
the homogeneity of functions hence that they commute with $E$.  So
they are homogeneous of degree $0$ for the $V$-filtration at $0$. On
the other hand, they are homogeneous of degree $1$ for the usual
filtration as any vector field. So if $u\in\gt(\kk)$,
$\gs_\gL\infun(u)=\gs(u)$.

On differential operators with constant coefficients, the
$V$-filtration at $\{0\}$ and the usual filtration coincide, hence we
have also $\gs^E\infun(P)=\gs(P)$ for these operators.

So, $Ch_\gL\infun(\CM_F)$ is contained in the set of
points where the symbols of the operators of $\gt(\kk)$ and of $F$
vanish that is in $(\kp\x\bN(\kp)) \cap\CC(\kp)$. 
\end{proof}

\begin{lemm} \label{lem:nilorbit}
For each nilpotent orbit $S$ of $\bN(\kp)$, there is a vector field
$\gh$ which is positive definite with respect to $S$ and such that
$\CM_F$ has $\gh$-weighted regular singularities along $T^*_S\kp$.
\end{lemm}

\begin{proof}
Let $S$ be one of these orbits, $r$ the codimension of $S$ and $X\in
S$. As in \cite[\S 3]{LEVASS2} (see also \cite[Part I, \S
5.6]{VARAD2}) we can choose a normal $\sld$-triple $(H,X,Y)$ in $\kp$
which generates a Lie algebra isomorphic to $\sld$ and acting on $\kp$
by the adjoint representation.  Then $\kp$ splits into a direct sum of
irreducible submodules of dimensions $\gl_i+1$ for $i=1\dots r$.
Moreover $\kp=\kp^Y\oplus[X,\kk]$, $\dim \kp^Y = r$ and we can select
a basis $(Y_1,\dots,Y_r)$ of $\kp^Y$ such that $[H,Y_i]=-\gl_iY_i$.
Let $V$ be a linear subspace of $\kk$ such that $\kk=V\oplus\kk^X$. If
$(b_1,\dots,b_{n-r})$ is a basis of $V$, the map $F:\C^n\to\kp$ given
by
$$F(x_1,\dots,x_{n-r},t_1,\dots,t_r)=
\exp(x_1b_1)\dots\exp(x_{n-r}b_{n-r}).(X+\sum t_iY_i)$$ is a local
isomorphism hence defines local coordinates $(x,t)$ of $\kp$ in a
neighborhood of $X$. In these coordinates, $S=\ensemble{(x,t)}{t=0}$,
$\kp^Y=\ensemble{(x,t)}{x=0}$, and the differential operators
$D_{x_1},\dots,D_{x_{n-r}}$ are in the ideal generated by $\gt(\kk)$
\cite[lemma 3.7]{SEKI}.

Let $E$ be the Euler vector field of the vector space $\kp$. A
standard calculation \cite[Part I,\S5.6]{VARAD2} shows that
$E(t_i)|_{x=0}=m_it_i$ with $m_i=\frac12\gl_i+1$. Moreover, if
$b_{n-r}=H$, we proved in \cite[lemma 3.4.1]{GL} that $E$ is equal to
$\gh+w$ where $\gh=\sum_{j=1}^{r}m_jt_jD_{t_j}$ and $w=1/2
D_{x_{n-r}}$. By definition, $\gh$ is positive definite with
respect to $S$.

Define a map $\gf:\kp\to V=\C^r$ by $\gf(x,t)=t$. Let
$\gh'=\sum m_jt_jD_{t_j}$ on $V$. The functions
$t_1,\dots,t_r$ satisfy $E(t_i)=\gh'(t_i)=m_it_i$ hence they are
homogeneous and the map $\gf$ is defined in a conic neighborhood of
$X$. This also shows that $E$ is invariant under $\gf$ and that
$\gh'=\gf_*(E)$.

The module $\CM_F$ is equal to $\Dp/\CI$ where $\CI$ is a coherent
ideal of $\Dp$ which contains the derivations
$D_{x_1},\dots,D_{x_{n-r}}$ hence $\CI$ is generated by
$D_{x_1},\dots,D_{x_{n-r}}$ and a finite set of differential operators
$Q_1(t,D_t),\dots,Q_N(t,D_t)$ depending only of $(t,D_t)$. (This
result is standard and also a special case of lemma \ref{lem:fin}).

The module $\CM_F$ satisfies the hypothesis of corollary
\ref{cor:transfer} hence $Ch^E\infun(\CM_F)$ is equal
to $Ch^{\gh}\infun(\CM_F)$ and by lemma \ref{lem:pp}
it is contained in $(\kp\x\bN(\kp)) \cap\CC(\kp)$.

Assume now that $T^*_S\kp$ is an irreducible component of the
characteristic variety $\CHMF$ and let $x^*$ be a generic point of
$T^*_S\kp$, that is a point which does not belong to other
irreducible components of $\CHMF$. We have
$T^*_S\kp\subset\CHMF\subset(\kp\x\bN(\kp))\cap\CC(\kp)$ and as they
have the same dimension, they are equal generically. So
$Ch^{\gh}\infun(\CM_F)=T^*_S\kp$ generically on
$T^*_S\kp$ and we are done.
\end{proof}

\begin{proof}[Proof of theorem \ref{thm:conj}]
We will argue by induction on the dimension of $\kg$ and first, we
reduce to the semi-simple case. Set $\kg_1=[\kg,\kg]$,
$\kk_1=\kk\cap\kg_1$, $\kp_1=\kp\cap\kg_1$, $\kz$ the center of $\kg$
and $\kp_0=\kz\cap\kp$. We have $\kp=\kp_0\oplus\kp_1$ and by lemma
\ref{lem:ss}, it suffices to prove the theorem for $\kp_1$. As
$\kz\cap\kk$ acts trivially we may assume that $\kg$ is semisimple.

Let $x$ be a non-nilpotent element of $\kp$. It  decomposes as $x=b+n$
where $b$ is non zero and semisimple, $n$ is nilpotent and
$[b,n]=0$. As $\kg$ is semisimple, $\kp^b$ is of dimension strictly
less than $\kp$, hence we may assume by the induction hypothesis that
the theorem is true for $\kp^b$. Lemma \ref{lem:closed} shows that
$\CM_F$ is regular in a neighborhood of $b$. As $\CM_F$ is constant on the orbits,
it is regular on the orbits whose closure contains $b$, in particular
at $x$.

We proved that $\CM_F$ is regular outside of the nilpotent cone. As the
nilpotent cone is a finite union of orbits, we will now argue by
descending induction on the dimension of these orbits. So let $x$ be a
nilpotent point of $\kp$, $Kx$ its orbit and assume that $\CM_F$ is
regular on $\kp-Kx$ in a neighborhood of $x$. Lemma \ref{lem:nilorbit}
shows that $\CM_F$ has $\gh$-weighted regular singularities along
$T^*_{Kx}\kp$ hence theorem \ref{thm:wregular} shows that $\CM_F$ is
regular at $x$.
\end{proof}

\backmatter


\begin{thebibliography}{10}

\bibitem{GL}
E.~Galina and Y.~Laurent, \emph{{D}-modules and characters of semi-simple {L}ie
  groups}, Pr{\'e}publications de l'Institut Fourier \textbf{570} (2002).

\bibitem{HC1}
Harish-Chandra, \emph{Invariant distributions on semi-simple {L}ie groups},
  Bull. Amer. Mat. Soc. \textbf{69} (1963), 117--123.

\bibitem{HC2}
Harish-Chandra, \emph{Invariant differential operators and distributions on a
  semi-simple {L}ie algebra}, Amer. J. Math. \textbf{86} (1964), 534--564.

\bibitem{HOTTA}
R.~Hotta and M.~Kashiwara, \emph{The invariant holonomic system on a semisimple
  lie algebra}, Inv. Math. \textbf{75} (1984), 327--358.

\bibitem{KVAN}
M.~Kashiwara, \emph{Vanishing cycles and holonomic systems of differential
  equations}, Lect. Notes in Math., vol. 1016, Springer, 1983, pp.~134--142.

\bibitem{KKREG}
M.~Kashiwara and T.~Kawa{\"\i}, \emph{On the holonomic systems of
  microdifferential equations {III}. systems with regular singularities}, Publ.
  RIMS, Kyoto Univ. \textbf{17} (1981), 813--979.

\bibitem{KOSR}
B.~Kostant and S.~Rallis, \emph{Orbits and representations associated with
  symmetric spaces}, Amer. J. Math \textbf{93} (1971), 753--809.

\bibitem{LAUM}
G.~Laumon, \emph{{$\mathcal D$}-modules filtr\'es}, Ast\'erisque, vol. 130,
  SMF, 1985, pp.~56--129.

\bibitem{THESE}
Y.~Laurent, \emph{Th\'eorie de la deuxi\`eme microlocalisation dans le domaine
  complexe}, Progress in Math., vol.~53, Birkh{\"a}user, 1985.

\bibitem{ENS}
Y.~Laurent, \emph{Polygone de {N}ewton et b-fonctions pour les modules
  microdiff\'erentiels}, Ann. Ec. Norm. Sup. 4e s\'erie \textbf{20} (1987),
  391--441.

\bibitem{LEVASS}
T.~Levasseur and J.T. Stafford, \emph{Invariant differential operators and a
  homomorphism of {H}arish-{C}handra}, Journal of the Americ. Math. Soc.
  \textbf{8} (1995), no.~2, 365--372.

\bibitem{LEVASS2}
T.~Levasseur and J.T. Stafford, \emph{Invariant differential operators on the tangent space of some
  symmetric spaces}, Ann. Inst. Fourier \textbf{49} (1999), no.~6, 1711--1741.

\bibitem{SABB}
C.~Sabbah, \emph{{$\mathcal D$}-modules et cycles \'evanescents}, G\'eom\'etrie
  r\'eelle, Travaux en cours, vol.~24, Hermann, 1987, pp.~53--98.

\bibitem{SKK}
M.~Sato, T.~Kawa{\"\i}, and M.~Kashiwara, \emph{Hyperfunctions and
  pseudo-differential equations}, Lect. Notes in Math., vol. 287, Springer,
  1980, pp.~265--529.

\bibitem{SEKI}
J.~Sekiguchi, \emph{Invariant spherical hyperfunctions on the tangent space of
  a symmetric space}, Advanced Studies in pure mathematics \textbf{6} (1985),
  83--126.

\bibitem{VARAD2}
V.S. Varadarajan, \emph{Harmonic analysis on real reductive groups}, Lect.
  Notes in Math., vol. 576, Springer, 1977.

\bibitem{WHIT}
H.~Whitney, \emph{Tangents to an analytic variety}, Annals of Math. \textbf{81}
  (1964), 496--549.

\end{thebibliography}
\enddocument

\end